\documentclass{amsart}
\usepackage{amssymb}
\usepackage{amsmath}
\usepackage{textcomp}
\makeindex

\newcommand{\ba}{\begin{array}}
\newcommand{\eea}{\end{eqnarray}}
\newcommand{\ea}{\end{array}}

\newtheorem{definition}{Definition}[section]
\newtheorem{theorem}[definition]{Theorem}

\usepackage[matrix,arrow,curve]{xy}
\usepackage[utf8x]{inputenc}
\usepackage{tikz}

\begin{document}
\title[On Existence Of non-compact Exact Lagrangian Cobordism]{On Existence Of Non-Compact Exact Lagrangian Cobordism}
\author[S. Mukherjee]{Sauvik Mukherjee}
\address{Presidency University, Kolkata, India.\\
e-mail:mukherjeesauvik@yahoo.com\\}

\begin{abstract} We prove an existence result for exact lagrangian cobordisms between closed legendrian.
\end{abstract}
\maketitle

\section{introduction} In \cite{Chekanov} and \cite{EkholmEtnyreSullivan}, counterexample has been provided to show that two closed legendrians which are formally legendrian isotopic need not be legendrian isotopic. However in \cite{Murphy} Murphy has proved that if the closed legendrians happens to be of very special kind namely, loose then they are legendrian isotopic.\\

On the other hand in \cite{Chantraine} Chantraine has proved that for two legendrian isopotic closed legendrians there must be an exact lagrangian cobordism joining them in the symplectization.\\

Combining these two result we get the following theorem.
\begin{theorem}(\cite{Murphy},\cite{Chantraine})
\label{MurphyChantraine}
If $L_1$ and $L_2$ be two closed loose legendrians in a contact manifold $(Y,\xi)$ and if $L_1$ and $L_2$ are formally legendrian isotopic then there exists an exact lagrangian cobordism joining them in the symplectization $Y\times \mathbb{R}$.
\end{theorem}

In this paper we prove that if $Y$ happens to open manifold and if we weaken the definition of exact lagrangian cobordism by removing the condition namely "the cobordism is compact after one removes the cylindrical ends", then we can drop the looseness hypothesis from \ref{MurphyChantraine}. So the main theorem of this paper is the following
\begin{theorem}
\label{Main}
If $L_1$ and $L_2$ be two closed legendrians (need not be loose) which are formally legendrian isotopic, in an open contact manifold $(Y^{2n+1},\xi),\ n\geq 2$ then there exists an exact lagrangian cobordism (definition differes from \cite{EkholmHondaKalman}) joining them in the symplectization $Y\times \mathbb{R}$.
\end{theorem}

\section{Preliminaries} In this section we define the preliminary notions and recall some known results that are going to be used in this paper. For definitions for legendrian isotopy and formal legendrian isotopy we refer to \cite{Murphy} and start by defining exact lagrangian cobordism. Our definition is little different from \cite{EkholmHondaKalman} as we shall not assume the cobordisms to be compact after removing the cylindrical ends.  
\begin{definition}
\label{Exactlagrangiancobordism}
First observe that the symplectization of contact manifold $(Y,\xi=ker\alpha)$ is $(Y\times \mathbb{R},d(e^t\alpha))$ is an exact symplectic manifold. A lagrangian submanifold $L\subset Y\times \mathbb{R}$ is called exact lagrangian if $e^t\alpha_{\mid L}$ is exact. An oriented exact lagrangian submanifolds $L\subset Y\times \mathbb{R}$ is called an exact lagrangian cobordism from a closed legendrians $L_1$ to a closed legendrian $L_2$ with cylindrical ends $\chi_1(L)$ and $\chi_2(L)$ if there exists $T>0$ such that
\[
\begin{array}{rclcl}
\chi_1(L)&=&L\cap (T,\infty)\times Y&=&(T,\infty)\times L_1\\
\chi_2(L)&=&L\cap (-\infty,-T)\times Y&=&(-\infty,-T)\times L_1\\
\end{array}
\]
and $f$ is constant on $\chi_i(L),\ i=1,2$ for $df=e^t\alpha_{\mid L}$.
\end{definition} 

Next we define the wrinkled legendrians as in \cite{Murphy} (no difference). Let us first recall wrinkled embeddings. Let $W:\mathbb{R}^n\to \mathbb{R}^{n+1}$ be a smooth proper map which is a topological embedding and which is a smooth embedding away from finitely many spheres $\{S^{n-1}_j\}$ and near the spheres $\{S^{n-1}_j\}$, $W$ is given by $W(u,v):=(v,x_1,z)$, where 
\[
\begin{array}{rcl}
x_1&=& u^3-3u(1-|v|^2)\\
z&=& (1/5)u^5-(2/3)u^3(1-|v|^2)+u(1-|v|^2)^2\\
\end{array}
\] 
where $u\in \mathbb{R}$ and $v\in \mathbb{R}^{n-1}$.
\begin{definition}
\label{Wrinkledlegendrian}
A topological embedding $f:L^n\to (Y^{2n+1},\xi)$ from a closed connected $L$ to a contact manifold $Y$ is called a wrinkled legendrian if $Im(df)\subset \xi$ everywhere and $df$ is of full rank outside finitely many $(n-2)$-spheres $\{S^{n-2}_j\}$, called legendrian wrinkles each one of which are contained in Darboux charts $U_j$, i.e, $S^{n-2}_j\subset U_j$, the composition of $f$ with the front projection on $L\cap U_j$ are wrinkled embeddings. For more details on wrinkled legendrians we refer to \cite{Murphy}.
\end{definition}

\begin{theorem}(\cite{Murphy})
\label{Murphy}
If $L_1$ and $L_2$ be closed legendrians (need not be loose) in a contact manifold $Y$ and if $L_1$ and $L_2$ are formally legendrian isotopic then there exists a isotopy of wrinkled legendrians joining $L_1$ with $L_2$.
\end{theorem}

\section{Legendrian Isotopy and Exact Lagrangian Cobordisms} In this section we outline the proof of a result from \cite{Chantraine}, namely 
\begin{theorem}(\cite{Chantraine})
\label{Chantraine}
If $L_0$ and $L_1$ be legendrians (need not be loose) in a contact manifold $(Y,\xi)$ which are legendrian isotopic then there exists an exact lagrangian cobordism joining them in the symplectization $Y\times \mathbb{R}$.
\end{theorem}
In order to sketch a proof of \ref{Chantraine} we need a theorem from \cite{Geiges}.
\begin{theorem}(\cite{Geiges})
\label{Geiges}
Let $L_t,\ t\in [0,1]$ be a legendrian isotopy of a closed legendrian $L$ in a contact manifold $(Y,\xi=ker\alpha)$. Then there exists a compactly supported contact isotopy $\psi_t:Y\to Y$ such that $\psi_t(L_0)=L_t$.
\end{theorem}
{\bf Sketch of the proof:} Define a one parameter family of vector fields along $L_t$ as $X_t:=\frac{d}{dt}(L_t)$. Now consider $\alpha(X_t)$ defined on $L_t$ and extend it to a tubular neighborhood of $L_t$ making it zero outside a slightly bigger tubular neighborhood. Call this function $H_t$. Now we just need to take the contact hamiltonian flow of of the contact hamiltonion vector field corresponding to the hamiltonian $H_t$. For more details we refer the readers to Theorem 2.41 of \cite{Geiges}.\\

 {\bf Sketch of the proof of \ref{Chantraine}:} Let $L_t$ be the legendrian isotopy between $L_0$ and $L_1$ and let $f_t$ be the compactly supported one parameter family of contactomorphisms as produced by \ref{Geiges}. In Proposition 2.2 in \cite{Chantraine} it has been proved that a contactomorphism of $Y$ can be lifted to a hamiltonian symplectomorphism of the symplectization $Y\times \mathbb{R}$. Let $H_t$ be the hamiltonian on $Y\times \mathbb{R}$ whose flow realizes the lift of $f_t$. Define $\tilde{H}:Y\times \mathbb{R}\times [0,1]\to \mathbb{R}$ as 
 \[
 \begin{array}{rcl}
 \tilde{H}(p,v,t)&=&H_t(p,v)\ for\ v>T\\
 \tilde{H}(p,v,t)&=&0\ for\ v<-T\\
 \end{array}
 \]
 Let $\phi_t$ be its hamiltonian flow. Now just consider $\phi_1(L_0\times \mathbb{R})$.\\
 
\section{Proof of \ref{Main}} It is clear from \ref{Murphy} and from the proof of \ref{Chantraine} that in order to prove \ref{Main} we need to understand how the function $H_t$ as in the proof of \ref{Geiges} behaves as a legendrian wrinkle appears and disappears. However one only needs to understand the phenomenon when a wrinkle appears, because when it disappears one only needs to replace $t$ by $-t$ for $t\in [0,1]$.\\

So consider the family of wrinkled legendrians whose front is given by $W_t:\mathbb{R}^2\to \mathbb{R}^3$ as $W_t(u,x_2)=(z,x_1,x_2)$, where 
\[
\begin{array}{rcl}
z&=&(1/5)u^5-(2/3)u^3(t-x_2^2)+u(t-x_2^2)^2\\
x_1&=&u^3-3u(t-x_2^2)\\
\end{array}
\]
Now set $y_1=(1/3)(u^2+x_2^2-t)$. Observe that $dz-y_1dx_1=[-(2/3)u^3x_2-2ux_2(t-x_2^2)]dx_2$. So we set $y_2=-(2/3)u^3x_2-2ux_2(t-x_2^2)$. For the higher dimensional case we just need to set $y_i=0\ for\ i>2$.\\

Now let us compute $X_t$. \[X_t=(-3u)\partial_{x_1}+(-1/3)\partial_{y_1}+(-2ux_2)\partial_{y_2}+(-(2/3)u^3+2u(t-x_2^2))\partial_{z}\] So for $\alpha=dz-y_dx_1-y_2dx_2$ we compute $\alpha(X_t)=(1/3)u^3+u(t-x_2^2)$.\\

In order to extend $\alpha(X_t)$ in a tubular neighborhood of $L_t$ we must remove $u$ from the formula of $\alpha(X_t)$. \\

Observe that $\alpha(X_t)=(x_1/3)+2u(t-x_2^2)$. First using the formula for $x_1$ and $y_2$ we get \[u=-\frac{(y_2+(2/3)x_1x_2)}{4x_2(t-x_2^2)}\] Using this formula for $u$ in the formula of $\alpha(X_t)$ we get $\alpha(X_t)=-(y_2/2x_2)$. So it gives the extension of $\alpha(X_t)$ only away from a tubular neighborhood of $L_t\cap \{x_2=0\}$ more precisely on $|x_2|\geq \varepsilon$.\\

Now solve the formula for $y_1$ for $u$ to get $u=\pm(3y_1-x_2^2+t)^{1/2}$. We substitute for $u=(3y_1-x_2^2+t)^{1/2}$ in the formula for $\alpha(X_t)$ on tubular neighborhood of $L_t\cap \{u\geq \delta\}$ and substitute $u=-(3y_1-x_2^2+t)^{1/2}$ for the formula for $\alpha(X_t)$ on a tubular neighborhood of $L_t\cap \{u\leq -\delta\}$.
 
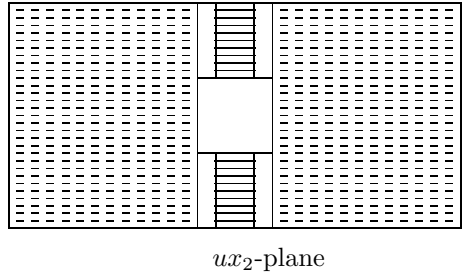
 \begin{center}
\begin{picture}(300,150)(-100,5)\setlength{\unitlength}{1cm}
\linethickness{.075mm}

\multiput(-1,1.5)(6,0){2}
{\line(0,1){3}}

\multiput(1.5,1.5)(1,0){2}
{\line(0,1){3}}

\multiput(1.75,1.5)(.5,0){2}
{\line(0,1){1}}

\multiput(1.75,3.5)(.5,0){2}
{\line(0,1){1}}

\multiput(-1,1.5)(0,3){2}
{\line(1,0){6}}

\multiput(1.5,2.5)(0,1){2}
{\line(1,0){1}}

\multiput(1.73,1.6)(.2,0){1}{\line(1,0){.55}}
\multiput(1.73,1.7)(.2,0){1}{\line(1,0){.55}}
\multiput(1.73,1.8)(.2,0){1}{\line(1,0){.55}}
\multiput(1.73,1.9)(.2,0){1}{\line(1,0){.55}}
\multiput(1.73,2)(.2,0){1}{\line(1,0){.55}}
\multiput(1.73,2.1)(.2,0){1}{\line(1,0){.55}}
\multiput(1.73,2.2)(.2,0){1}{\line(1,0){.55}}
\multiput(1.73,2.3)(.2,0){1}{\line(1,0){.55}}
\multiput(1.73,2.4)(.2,0){1}{\line(1,0){.55}}

\multiput(1.73,3.6)(.2,0){1}{\line(1,0){.55}}
\multiput(1.73,3.7)(.2,0){1}{\line(1,0){.55}}
\multiput(1.73,3.8)(.2,0){1}{\line(1,0){.55}}
\multiput(1.73,3.9)(.2,0){1}{\line(1,0){.55}}
\multiput(1.73,4)(.2,0){1}{\line(1,0){.55}}
\multiput(1.73,4.1)(.2,0){1}{\line(1,0){.55}}
\multiput(1.73,4.2)(.2,0){1}{\line(1,0){.55}}
\multiput(1.73,4.3)(.2,0){1}{\line(1,0){.55}}
\multiput(1.73,4.4)(.2,0){1}{\line(1,0){.55}}

\multiput(2.6,1.6)(.2,0){12}{\line(1,0){.09}}
\multiput(2.6,1.7)(.2,0){12}{\line(1,0){.09}}
\multiput(2.6,1.8)(.2,0){12}{\line(1,0){.09}}
\multiput(2.6,1.9)(.2,0){12}{\line(1,0){.09}}
\multiput(2.6,2)(.2,0){12}{\line(1,0){.09}}
\multiput(2.6,2.1)(.2,0){12}{\line(1,0){.09}}
\multiput(2.6,2.2)(.2,0){12}{\line(1,0){.09}}
\multiput(2.6,2.3)(.2,0){12}{\line(1,0){.09}}
\multiput(2.6,2.4)(.2,0){12}{\line(1,0){.09}}
\multiput(2.6,2.5)(.2,0){12}{\line(1,0){.09}}
\multiput(2.6,2.6)(.2,0){12}{\line(1,0){.09}}
\multiput(2.6,2.7)(.2,0){12}{\line(1,0){.09}}
\multiput(2.6,2.8)(.2,0){12}{\line(1,0){.09}}
\multiput(2.6,2.9)(.2,0){12}{\line(1,0){.09}}
\multiput(2.6,3)(.2,0){12}{\line(1,0){.09}}
\multiput(2.6,3.1)(.2,0){12}{\line(1,0){.09}}
\multiput(2.6,3.2)(.2,0){12}{\line(1,0){.09}}
\multiput(2.6,3.3)(.2,0){12}{\line(1,0){.09}}
\multiput(2.6,3.4)(.2,0){12}{\line(1,0){.09}}
\multiput(2.6,3.5)(.2,0){12}{\line(1,0){.09}}
\multiput(2.6,3.6)(.2,0){12}{\line(1,0){.09}}
\multiput(2.6,3.7)(.2,0){12}{\line(1,0){.09}}
\multiput(2.6,3.8)(.2,0){12}{\line(1,0){.09}}
\multiput(2.6,3.9)(.2,0){12}{\line(1,0){.09}}
\multiput(2.6,4)(.2,0){12}{\line(1,0){.09}}
\multiput(2.6,4.1)(.2,0){12}{\line(1,0){.09}}
\multiput(2.6,4.2)(.2,0){12}{\line(1,0){.09}}
\multiput(2.6,4.3)(.2,0){12}{\line(1,0){.09}}
\multiput(2.6,4.4)(.2,0){12}{\line(1,0){.09}}

\multiput(-.9,1.6)(.2,0){12}{\line(1,0){.09}}
\multiput(-.9,1.7)(.2,0){12}{\line(1,0){.09}}
\multiput(-.9,1.8)(.2,0){12}{\line(1,0){.09}}
\multiput(-.9,1.9)(.2,0){12}{\line(1,0){.09}}
\multiput(-.9,2)(.2,0){12}{\line(1,0){.09}}
\multiput(-.9,2.1)(.2,0){12}{\line(1,0){.09}}
\multiput(-.9,2.2)(.2,0){12}{\line(1,0){.09}}
\multiput(-.9,2.3)(.2,0){12}{\line(1,0){.09}}
\multiput(-.9,2.4)(.2,0){12}{\line(1,0){.09}}
\multiput(-.9,2.5)(.2,0){12}{\line(1,0){.09}}
\multiput(-.9,2.6)(.2,0){12}{\line(1,0){.09}}
\multiput(-.9,2.7)(.2,0){12}{\line(1,0){.09}}
\multiput(-.9,2.8)(.2,0){12}{\line(1,0){.09}}
\multiput(-.9,2.9)(.2,0){12}{\line(1,0){.09}}
\multiput(-.9,3)(.2,0){12}{\line(1,0){.09}}
\multiput(-.9,3.1)(.2,0){12}{\line(1,0){.09}}
\multiput(-.9,3.2)(.2,0){12}{\line(1,0){.09}}
\multiput(-.9,3.3)(.2,0){12}{\line(1,0){.09}}
\multiput(-.9,3.4)(.2,0){12}{\line(1,0){.09}}
\multiput(-.9,3.5)(.2,0){12}{\line(1,0){.09}}
\multiput(-.9,3.6)(.2,0){12}{\line(1,0){.09}}
\multiput(-.9,3.7)(.2,0){12}{\line(1,0){.09}}
\multiput(-.9,3.8)(.2,0){12}{\line(1,0){.09}}
\multiput(-.9,3.9)(.2,0){12}{\line(1,0){.09}}
\multiput(-.9,4)(.2,0){12}{\line(1,0){.09}}
\multiput(-.9,4.1)(.2,0){12}{\line(1,0){.09}}
\multiput(-.9,4.2)(.2,0){12}{\line(1,0){.09}}
\multiput(-.9,4.3)(.2,0){12}{\line(1,0){.09}}
\multiput(-.9,4.4)(.2,0){12}{\line(1,0){.09}}

\put(1.7,1){$ux_2$-plane}

\end{picture}\end{center}

In the above picture the big shaded rectangles represents $|u|\geq\delta$ and the central vertical striped strips corresponds to $\{|x_2|\geq \varepsilon\}\cap \{|u|\leq \delta/2\}$. The central non-shaded rectangle represents $\{|u|<\delta\}\cap \{|x_2|<\varepsilon\}$. Now on vertical non-shaded strips we linearly homotope from one definition of $u$ to the other which agree along the legendrian $L_t$ and disagree in the normal direction of the tubular neighborhood of $L_t$.\\

We shall actually replace the non-shaded rectangle by a circular one namely by $\{u^2+x_2^2<\varepsilon\}$.\\

Now we shall follow a technique explained in \cite{Ekholm} which allows us to push the non-shaded central rectangle (disc) to infinity in $Y$ along a family of isotropic path which begins on the non-shaded rectangle (disc) and goes to infinity in $Y$.\\

Let us rename $L_1$ and $L_2$ in \ref{Main} by $L_{-T}$ and $L_T$ respectively and call the wrinkle legendrian isotopy between $L_{-T}$ and $L_T$ by $L_t,\ t\in [-T,T]$. Define $g:[-3T,3T]\to \mathbb{R}$ as 
\[
\begin{array}{rcl}
g(t)&=&\infty\ for\ t\in [-T,T]\\
g(t)&=&0\ for\ t\in [-3T,-2T]\cup [2T,3T]\\
\end{array}
\]   

Consider the non-shaded central rectangle (disc) and its image for different values of $t\in [-3T,3T]$. Any point of this set we denote by $p_t$.\\  

Now define an isotropic path depending on $p_t$, $\gamma_{p_t}:[0,g(t)]\to Y$ as follows. The $z$-component of $\gamma_{p_t}(\tau)$ is $\tau y_1^t+z^t$, where $y_1^t$ is the $y_1$ component of $p_t$, $z^t$ is the $z$ component of $p_t$ and the $x_1$ component of $\gamma_{p_t}(\tau)$ is $\tau+x_1^t$, where $x_1^t$ is the $x_1$ component of $p_t$ and all other component of $\gamma_{p_t}(\tau)$ is constant and equal to the component of $p_t$. Clearly it is an isotropic path. Observe that if we replace $\tau$ by $-\tau$ for the components of $\gamma_{p_t}(\tau)$ we shall get an isotropic path extended in the reverse direction. So we can choose the path so that it does not intersect the legendrian for $\tau >0$.\\

The lagrangian projection of this path is an isotropic path in $\mathbb{R}^4$ ($\mathbb{R}^{2n}$). As in section 3.2.3 of \cite{Ekholm} we can choose a tubular neighborhood of this projected path of the form $(-\varepsilon,g(t)+\varepsilon)\times (-\varepsilon,\varepsilon)$ and for higher dimensional case of the form ($(-\varepsilon,g(t)+\varepsilon)\times (-\varepsilon,\varepsilon)\times B_{\varepsilon}^{2n-2}$), where $B_{\varepsilon}^{2n-2}$ is the standard symplectic ball of diameter $2\varepsilon$. Define as in \cite{Ekholm} $G(\tau,s,a,b)=\psi(s)f(\tau)\phi(a,b)$, where we have 
\begin{enumerate}
\item $\phi(a,b)=1,\ for\ |(a,b)|\leq (1/4)\varepsilon,\ \phi(a,b)=0,\ for\ |(a,b)|\geq (3/4)\varepsilon$ and $|d\phi(a,b)|=\mathcal{O}(\varepsilon^{-1})$\\
\item $f(\tau)=1,\ for\ \tau\in [0,g(t)],\ f(\tau)=0,\ for\ \tau \notin [(-\varepsilon/2),g(t)+(-\varepsilon/2)],\ and\ |f'(\tau)|=\mathcal{O}(\varepsilon^{-1})$\\
\item $\psi(s)=\mathcal{O}(\varepsilon^{-1}),\ \psi'(s)=1\ for |s|\leq \varepsilon/4\ and\  \psi(s)=\psi'(s)=0,\ |s|>(3/4)\varepsilon$
\end{enumerate} 

Take the contact lift of the time $g(t)$ flow of $X_G$ and call it $\Psi_t$. Now $\alpha(X_t)$ is changed to $\alpha(d\Psi_t(X_t))$. But $\Psi_t$ are contactomorphisms and hence $\alpha(d\Psi_t(X_t))=\Psi_t^*\alpha(X_t)=\lambda(p_t)\alpha(X_t)$. Now $\lambda(p_t)$ is infinity for $t\in [-T,T]$. Observe that if a wrinkle appears for $t=t_0$ and disappears for $t=t_1$ then $[t_0,t_1]\subset [-T,T]$. \\

So in view of the above argument we shall prove \ref{Main} if we explain the situation of nested wrinkles. But observe that the size of the non-shaded central rectangle (disc) which we have pushed to infinity in $Y$ could be taken as small as we want. So this region for external wrinkles are chosen so small that the regions corresponding to the inner most wrinkle contains it. The following picture explains the situation and completes the proof of \ref{Main} (Obviously one needs to replace $T$ by $3T$ in the definition of exact lagrangian cobordisms \ref{Exactlagrangiancobordism} ).

 \begin{center}
\begin{picture}(300,150)(-100,5)\setlength{\unitlength}{1cm}
\linethickness{.075mm}

\multiput(-1,1.5)(6,0){2}
{\line(0,1){3}}

\multiput(-1,1.5)(0,3){2}
{\line(1,0){6}}

\qbezier(-.5,3)(2,1)(4.5,3)
\qbezier(-.5,3)(2,5)(4.5,3)

\qbezier(2,2.2)(1,2.8)(2,3.8)
\qbezier(2,2.2)(3,2.8)(2,3.8)

\qbezier(1,2.5)(1.5,2.8)(1,3.3)
\qbezier(1,2.5)(.5,2.8)(1,3.3)

\qbezier(3,2.5)(3.5,2.8)(3,3.3)
\qbezier(3,2.5)(2.5,2.8)(3,3.3)

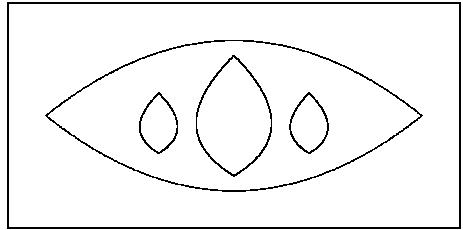
\put(0,1){Singular locus of nested wrinkles}

\end{picture}\end{center}

In the above picture the spherical curves represents the singular locus of nested wrinkles. Now the region that we pushed to infinity for the outer most wrinkle can be chosen so small that the region for the central inner wrinkle contains it. However for the regions corresponding to the left and right inner wrinkles needs to be considered separately. So we only need to consider the inner most wrinkles and not the outer ones.

\end{document}